# A Dynamized Power Flow Method based on Differential Transformation

Yang Liu, *Student Member*, *IEEE*, Kai Sun, *Senior Member*, *IEEE*, Jiaojiao Dong, *Member*, *IEEE*

*Abstract*— This paper proposes a novel method for solving and tracing power flow solutions with changes of a loading parameter. Different from the conventional continuation power flow method, which repeatedly solves static AC power flow equations, the proposed method extends the power flow model into a fictitious dynamic system by adding a differential equation on the loading parameter. As a result, the original solution curve tracing problem is converted to solving the time domain trajectories of the reformulated dynamic system. A non-iterative algorithm based on differential transformation is proposed to analytically solve the aforementioned dynamized model in form of power series of time. This paper proves that the nonlinear power flow equations in the time domain are converted to formally linear equations in the domain of the power series order after the differential transformation, thus avoiding numerical iterations. Case studies on several test systems including a 2383-bus system show the merits of the proposed method.

## I. INTRODUCTION

TRACING solution curves of power flow equations with the changes of a trending parameter such as the loading level is usually computation-intensive task in power system operations and planning to prevent steady-state voltage instability and other insecurities [1]-[4]. An example is computation of the power-voltage (P-V) curves for critical buses with the increase of load. Traditionally, the continuation power flow (CPF) method [5]-[9] is widely used to solve P-V curves, which adopts a prediction-correction scheme to identify a series of power flow solutions along the solution curve where each prediction step gives an initial guess and the following corrector step performs numerical iterations to find the converged solution [10]-[12]. However, the CPF method may suffer from huge computation burdens since it requires solving nonlinear AC power flow equations for multiple times using numerical iteration methods [13]-[14]. Moreover, the computation burden of the CPF method can further grow and become unacceptable with modern power grids being integrated with high renewable energy resources and demand response, where such a solution process with numerical iterations needs to be repeated many times for multiple contingencies.

In the literature, some techniques are proposed to reduce the computation burden of the CPF method [15]-[17], falling into two categories. The first category of techniques aim to design a more effective predictor than a standard tangent predictor [10]-[12] as adopted by many commercial CPF programs [13]-[14]. For example, paper [15] proposes three types of nonlinear predictors to predict a new solution from more than one previous solutions using interpolation and polynomial approximation including the Lagrange's polynomial interpolation formula, Newton's forward and backward divided difference formula, and cubic spline interpolation method. Besides, a secant predictor is used in [16] and a holomorphic embedding-based predictor is proposed in [17]. Methods in the first category are able to generate a better initial guess for the Newton Raphson (NR) method so as to reduce the total number of iterations; however, they still require many numerical iterations for solving nonlinear power flow equations. The second category of techniques focus on more efficient correctors than the standard NR method-based corrector. For example, authors in [15] propose a hybrid corrector allowing the switches between a NR method (taking its merit of robustness) and a fast decoupled NR method (taking its merit of fast speed) until a pre-defined maximum total number of iterations. Methods in this category are mainly used to improve the convergence performance of numerical iterations, but still require solving nonlinear AC power flow equations repeatedly. Overall, a major bottleneck for the above two categories of methods lies in their inherent solution mechanism that the power-flow equations are essentially solved as algebraic equations in an iterative manner.

To more efficiently trace solution curves of power flow equations, this paper proposes a novel dynamized power flow (DPF) method that extends the power flow model into a fictitious dynamic system, called a "dynamized" power flow model, by adding a differential equation about a fictitious time, and then solve the complete time-domain trajectory of the dynamic system instead of repeatedly solving power flow equations for a series of conditions. A differential transformation (DT) method, which is proved effective for solving power system transient stability simulation in our recent works [18]-[20], is applied to solve the dynamized model, named as dynamized power flow (DPF) method. This paper proves that the nonlinear AC power flow equations are converted to formally linear equations after DT, and further

This work was supported in part by the ERC Program of the NSF and DOE under NSF Grant EEC-1041877 and in part by NSF Grant ECCS-1610025.
Y. Liu, K. Sun and J. Dong are with the Department of EECS, University of Tennessee, Knoxville, TN 37996 USA (e-mail: yliu161@vols.utk.edu, kaisun@utk.edu, jdong7@utk.edu).

designs an efficient algorithm to solve the time domain trajectory without numerical iterations. Case studies on several test systems including a 2383-bus system demonstrate the accuracy, computational complexity and time performance of the proposed approach compared with a CPF solver.

The rest of the paper is organized as follows. Section II gives the problem description. Section III presents the proposed method. Section IV is case study and Section V draws conclusions.

## II. PROBLEM STATEMENT

The conventional power flow equations are given in (1a) where $\overline{\mathbf{S}}$ is a vector of the complex power injections, $\overline{\mathbf{V}}$ is a vector of bus voltage phasor, and $\mathbf{Y_{bus}}$ is the bus admittance matrix. By adding the product of a loading parameter $\lambda$ and a constant vector $\overline{\mathbf{b}}$ to the left-hand side, a general continuum of power flow equations is formulated in (1b).

$$\begin{array}{ll} \overline{\mathbf{S}} = \overline{\mathbf{V}}(\mathbf{Y_{bus}}\overline{\mathbf{V}})^* & \text{(a)} \\ \overline{\mathbf{S}} + \lambda\overline{\mathbf{b}} = \overline{\mathbf{V}}(\mathbf{Y_{bus}}\overline{\mathbf{V}})^* & \text{(b)} \end{array} \quad (1)$$

Note that the vector $\overline{\mathbf{b}}$ is defined to reflect an arbitrarily direction of load changes, for example, uniform increases of all generation and load, or increases of generation and load at certain buses or zones. Meanwhile, practical operating constraints such as the reactive power limit of generators can be considered during the load change.

Equation (1b) is further written as the general form in (2) where $\boldsymbol{g}$ is a nonlinear vector field; $\boldsymbol{y}$ is the bus voltage vector under rectangular coordinates defined as $\boldsymbol{y}=[\boldsymbol{e}^{\mathrm{T}},\boldsymbol{f}^{\mathrm{T}}]^{\mathrm{T}}$, where $\boldsymbol{e}=[e_1,\ldots,e_N]^{\mathrm{T}}$ and $\boldsymbol{f}=[f_1,\ldots,f_N]^{\mathrm{T}}$ are respectively the real and imaginary parts of the bus voltage phasor; $N$ is the total number of buses; $\lambda$ is the loading parameter.

$$0 = \boldsymbol{g}(\boldsymbol{y},\lambda) \quad (2)$$

The goal is to determine how power flow solution $\boldsymbol{y}$ changes with loading parameter $\lambda$, shown in (3). After (3) is obtained, the other system variables (such as voltage magnitude and power injections) are easily calculated to draw P-V curves.

$$\boldsymbol{y} = \boldsymbol{y}(\lambda) \quad (3)$$

Generally, analytical expression of (3) is unavailable due to the nonlinearity of $\boldsymbol{g}$ in (2). Therefore, a prediction-correction scheme and numerical iterations are needed in conventional CPF method.

## III. PROPOSED DYNAMIZED POWER FLOW METHOD

### A. Introduction of the Differential Transformation

A smooth nonlinear function of time $x(t)$ can be approximated by a $K^{\text{-th}}$ order polynomial function of time as shown in (4), where $X(k)$ is the $k^{\text{th}}$ order power series coefficient and can be calculated by (5).

$$x(t) = \sum_{k=0}^{K} X(k)t^k \quad (4)$$

$$X(k) = \frac{1}{k!}\left[\frac{d^k x(t)}{dt^k}\right]_{t=0} \quad (5)$$

Generally, these power series coefficients are calculated in a recursive manner from $k=0$ to $k=K$, and many mathematical methods can be used such as the Adomain decomposition method [21]-[22] and the power series-based method in [23]-[24]. However, the applications of the above methods are limited by their huge computational burdens in deriving power series coefficient $X(k)$ using the complicated high order derivative operations.

As an emerging mathematical tool, DT [25]-[28] considers power series coefficient $X(k)$ as a transformation of $x(t)$ at the $k^{\text{th}}$ order as shown by (6). When multiple functions like $x(t)$ are to be calculated and analyzed, their high order power series coefficients can directly be operated and calculated based on transformation rules introduced by DT. Thus, there is no need to calculate complicated high order derivatives of each function.

$$x(t) \rightarrow X(k) \quad (6)$$

Our recent paper [18]-[19] introduces DT to the power system field to effectively solve power system nonlinear differential-algebraic equations (DAEs) for transient stability simulation. New transformation rules for nonlinear functions in power system models are proved in [18]-[19]. Five rules are given in (7) and will be utilized in this paper. Here, $X(k)$ and $Y(k)$ are DTs of functions $x(t)$ and $y(t)$, $c$ is a constant, and $\delta$ is the Kronecker delta function:

i) $x(t) \pm y(t) \rightarrow X(k) \pm Y(k)$;   ii) $cx(t) \rightarrow cX(k)$

iii) $x(t)y(t) \rightarrow X(k) \otimes Y(k) \triangleq \sum_{m=0}^{k} X(m)Y(k-m)$   (7)

iv) $\dfrac{dx(t)}{dt} \rightarrow (k+1)X(k+1)$; v) $c \rightarrow c\delta(k)=\begin{cases}1, k=0\\0, k\neq 0\end{cases}$

### B. Conceptual Description of the Proposed Method

The proposed method has following four steps, where each step is first briefed below and then described in detail in Section III-C to Section III-F respectively.

First, the algebraic equation (2) is extended to a set of DAEs by introducing a fictitious time $t$ and adding two new equations, i.e., (8a) and (8b). Differential equation (8a) is a dynamic system to trace the changes of system variables such as power or voltages, where $x(t)$ is a state variable and $f(\cdot)$ is a vector field. Algebraic equation (8b) is an ancillary equation that builds the relationship between the newly introduced variable $x(t)$ and the original variables $\boldsymbol{y}(t)$ and $\lambda(t)$. Note that the ancillary function $h$ may not be needed if $x(t)$ is selected from one of the variables in $\boldsymbol{y}(t)$ and $\lambda(t)$. The details of designing (8a) and (8b) are described in Section III-C.

$$\begin{array}{ll} \dot{x}(t) = f(x(t),\boldsymbol{y}(t),\lambda(t)) & \text{(a)} \\ 0 = h(x(t),\boldsymbol{y}(t),\lambda(t)) & \text{(b)} \\ 0 = \boldsymbol{g}(\boldsymbol{y}(t),\lambda(t)), i.e., \; Equ. \; (2) & \text{(c)} \end{array} \quad (8)$$

Second, the DTs of (8a)-(8c) are derived in (9a)-(9c) respectively, using the transformation rules in (7). Specifically, the nonlinear power flow equation (8c) is converted to a new set of equations (9c) that couples the power series coefficients of $\boldsymbol{y}(t)$ and $\lambda(t)$ in all orders, i.e., $\boldsymbol{Y}(0)\ldots \boldsymbol{Y}(k), \Lambda(0)\ldots \Lambda(k)$.

$$(k+1)X(k+1) = F(X(0:k), \boldsymbol{Y}(0:k), \Lambda(0:k)) \quad (a)$$
$$0 = H(X(0:k), \boldsymbol{Y}(0:k), \Lambda(0:k)) \quad (b) \quad (9)$$
$$0 = \boldsymbol{G}(\boldsymbol{Y}(0:k), \Lambda(0:k)) \quad (c)$$

Third, we prove that both (9c) and (9b) satisfy formally linear equations about the $k^{th}$ order power series coefficients $\boldsymbol{Y}(k)$ and $\Lambda(k)$, as shown in (10a) and (10b), respectively, where $A$ matrices are functions of $\boldsymbol{Y}(0)$ and $\Lambda(0)$ and $B$ matrices are functions of $\boldsymbol{Y}(0:k-1)$ and $\Lambda(0:k-1)$. As a result, $\boldsymbol{Y}(k)$, i.e., the $k^{th}$ order power series coefficient of bus voltage vector, is analytical solved from $\boldsymbol{Y}(0:k-1)$ and $\Lambda(0:k-1)$, either by (11) or by (12), depending on if the ancillary function $h$ is needed when designing the differential equation in (8).

$$0 = \boldsymbol{A}_{gy}\boldsymbol{Y}(k) + \boldsymbol{A}_{g\lambda}\Lambda(k) + \boldsymbol{B}_g \quad (a)$$
$$0 = \boldsymbol{A}_{hy}\boldsymbol{Y}(k) + \boldsymbol{A}_{h\lambda}\Lambda(k) + \boldsymbol{B}_h \quad (b) \quad (10)$$

$$\boldsymbol{Y}(k) = -\boldsymbol{A}_{gy}^{-1}(\boldsymbol{A}_{g\lambda}\Lambda(k) + \boldsymbol{B}_g) \quad (11)$$

$$\begin{bmatrix} \boldsymbol{Y}(k) \\ \Lambda(k) \end{bmatrix} = -\begin{bmatrix} \boldsymbol{A}_{gy} & \boldsymbol{A}_{g\lambda} \\ \boldsymbol{A}_{hy} & \boldsymbol{A}_{h\lambda} \end{bmatrix}^{-1} \begin{bmatrix} \boldsymbol{B}_g \\ \boldsymbol{B}_h \end{bmatrix} \quad (12)$$

Finally, we design a non-iterative algorithm based on (9a) and (11) or (12) to solve power series coefficients $X(k)$, $\boldsymbol{Y}(k)$ and $\Lambda(k)$ from $k=0$ to any order $K$ in a recursively manner, and approximate variables $x(t)$, $\boldsymbol{y}(t)$ and $\lambda(t)$ as power series of time, shown in (13). After $\boldsymbol{y}(t)$ and $\lambda(t)$ are solved, the solution curves of power flow equations are directly obtained, as illustrated in Section IV-A.

$$x(t) = X(0) + X(1)t + X(2)t^2 + ...X(K)t^K$$
$$\boldsymbol{y}(t) = \boldsymbol{Y}(0) + \boldsymbol{Y}(1)t + \boldsymbol{Y}(2)t^2 + ...\boldsymbol{Y}(K)t^K \quad (13)$$
$$\lambda(t) = \Lambda(0) + \Lambda(1)t + \Lambda(2)t^2 + ...\Lambda(K)t^K$$

Among the above four steps, only the last step needs to be performed online, while the first three steps can be conducted in the offline stage because they are mainly used to derive expressions of matrices $A$ and $B$ in (10) and function $F$ in (9a), which is a one-time effort.

**Remarks:** there are two important observations: 1) from (10a) that the nonlinear power flow equation (2) about $\boldsymbol{y}(t)$ are converted to a formally linear equation about power series coefficients $\boldsymbol{Y}(k)$ after DT; 2) coefficients on bus voltages are explicitly solved by (11) or (12) and then used to calculate bus voltages by (13) in a straightforward manner, which is different from using a conventional power flow solver to calculate bus voltages by numerical iterations. The proposed DT based method for solving solution curves of power flow equations differentiates itself from the traditional continuation power flow method that relies on iterative numerical methods such as the family of Newton Raphson methods.

### C. Step 1: Dynamization of the Power Flow Equation

Two formulations of (8) are proposed to dynamize the power flow equation (2), shown in (14) and (15) respectively, where $C_1$ and $C_2$ are constants and $v_l(t)$ is the voltage magnitude of a load bus $l$. In (14), there is no ancillary equation (8b) because the selected state variable $\lambda(t)$ has existed in (2). In (15), the ancillary equation gives the relationship between bus voltage magnitude and the rectangular coordinate components.

**Formulation 1:**
$$\dot{\lambda}(t) = C_1$$
$$0 = \boldsymbol{g}(\boldsymbol{y}(t), \lambda(t)), i.e., Equ.(2) \quad (14)$$

**Formulation 2:**
$$\dot{v}_l(t) = C_2$$
$$0 = v_l(t)^2 - e_l(t)^2 - f_l(t)^2 \quad (15)$$
$$0 = \boldsymbol{g}(\boldsymbol{y}(t), \lambda(t)), i.e., Equ.(2)$$

For Formulation 1, its purpose is to characterize how the power changes with time, i.e., the power increases with time when $C_1>0$ and decreases with time when $C_1<0$. It can be used to trace curve segment in various shapes, either monotonically or non-monotonically, such as the curves (a), (b) and (c) in Fig. 1. For Formulation 2, its purpose is to characterize how the voltage magnitude changes with time, i.e., the voltage magnitude increases with time when $C_2>0$ and decreases with time when $C_2<0$. It can also be used to trace either monotonical or non-monotonical curve segments such as (a), (b) and (d) in Fig. 1.

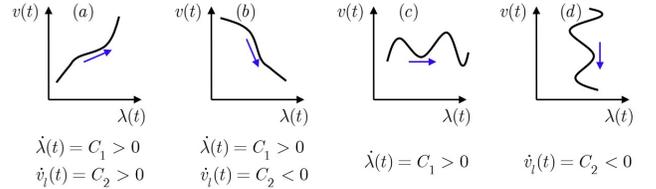

Fig. 1. Illustration of the two dynamized formulations for tracing curve segments of power flow equations

The above two formulations can be flexibly used to trace the full solution curve of a power flow equation. For example, the high voltage solutions in a P-V curve can be traced by Formulation 1 with $C_1>0$, the low voltage solutions can be traced by Formulation 1 with $C_1<0$, and the solution curves near the nose point can be traced by Formulation 2 with $C_2<0$.

### D. Step 2: Deriving DTs of the Dynamized Power Flow Model

#### 1) DTs of Nonlinear Power Flow Equation

The nonlinear power flow equation (2) is written into (16)-(19) under rectangular coordinates, where $\Omega_{PQ}$, $\Omega_{PV}$, $\Omega_{REF}$ are the set of PQ buses, PV buses and reference bus respectively, $p$ and $q$ are active and reactive power, $e$ and $f$ are the real and imaginary parts of bus voltages, $g$ and $b$ are real and imaginary parts of the admittance, $v$ is the voltage magnitude, superscript $sp$ means the value is specified, subscript $i$ and $j$ are the index of buses.

$$p_i^{sp} = g_p(\boldsymbol{y}, \lambda)$$
$$= -\lambda \Delta p_i + \sum_{j=1}^{N} g_{ij}\left(e_i e_j + f_i f_j\right) + \sum_{j=1}^{N} b_{ij}\left(f_i e_j - e_i f_j\right) \quad (16)$$
$$\text{if } i \in \Omega_{PQ} \cup \Omega_{PV}$$

$$q_i^{sp} = g_q(\boldsymbol{y}, \lambda)$$

$$= -\lambda \Delta q_i - \sum_{j=1}^{N} b_{ij} \left( e_i e_j + f_i f_j \right) + \sum_{j=1}^{N} g_{ij} \left( f_i e_j - e_i f_j \right) \quad (17)$$
$$\text{if } i \in \Omega_{PQ}$$

$$(v_i^{sp})^2 = g_v(\boldsymbol{y}) = e_i^2 + f_i^2, \text{if } i \in \Omega_{PV} \quad (18)$$

$$e_i^{sp} = g_e(\boldsymbol{y}) = e_i, \text{if } i \in \Omega_{REF} \quad (19)$$
$$f_i^{sp} = g_f(\boldsymbol{y}) = f_i, \text{if } i \in \Omega_{REF}$$

The DTs of (16)-(19) are derived in (20)-(23), respectively.

$$p_i^{sp}\delta(k) = G_p(\boldsymbol{Y}, \boldsymbol{\Lambda})$$

$$= -\Delta p_i \Lambda(k) + \sum_{j=1}^{N} g_{ij} \left( E_i(k) \otimes E_j(k) + F_i(k) \otimes F_j(k) \right)$$
$$+ \sum_{j=1}^{N} b_{ij} \left( F_i(k) \otimes E_j(k) - E_i(k) \otimes F_j(k) \right)$$
$$\text{if } i \in \Omega_{PQ} \cup \Omega_{PV}$$
$$(20)$$

$$q_i^{sp}\delta(k) = G_q(\boldsymbol{Y}, \boldsymbol{\Lambda})$$

$$= -\Delta q_i \Lambda(k) - \sum_{j=1}^{N} b_{ij} \left( E_i(k) \otimes E_j(k) + F_i(k) \otimes F_j(k) \right)$$
$$+ \sum_{j=1}^{N} g_{ij} \left( F_i(k) \otimes E_j(k) - E_i(k) \otimes F_j(k) \right)$$
$$\text{if } i \in \Omega_{PQ}$$
$$(21)$$

$$(v_i^{sp})^2 \delta(k) = G_v(\boldsymbol{Y})$$
$$= E_i(k) \otimes E_i(k) + F_i(k) \otimes F_i(k), \text{if } i \in \Omega_{PV} \quad (22)$$

$$e_i^{sp}\delta(k) = G_e(\boldsymbol{Y}) = E_i(k), \text{if } i \in \Omega_{REF}$$
$$f_i^{sp}\delta(k) = G_f(\boldsymbol{Y}) = F_i(k), \text{if } i \in \Omega_{REF} \quad (23)$$

For details, the derivation of (20) is elaborated below as an example. The remaining equations (21)-(23) are obtained in a similar procedure.

The left-hand-side (LHS) and the first term in the right-hand-side (RHS) of (20) are obtained by applying the rule (7-i), (7-ii) and (7-v) to the corresponding terms of (16). Note that $p_i^{sp}$ and $\Delta p_i$ are constants and $\lambda = \lambda(t)$ is a variable, therefore their transformations are: $p_i^{sp} \rightarrow p_i^{sp}\delta(k)$ and $\Delta p_i \lambda \rightarrow \Delta p_i \Lambda(k)$.

The remaining terms in the RHS of (20) are obtained from the corresponding terms of (16) by following steps:

First, apply the rule (7-iii):

$$e_i e_j \rightarrow E_i(k) \otimes E_j(k) \quad f_i f_j \rightarrow F_i(k) \otimes F_j(k)$$
$$f_i e_j \rightarrow F_i(k) \otimes E_j(k) \quad e_i f_j \rightarrow E_i(k) \otimes F_j(k)$$

Then, apply the rule (7-i) and (7-ii):

$$g_{ij} \left( e_i e_j + f_i f_j \right) \rightarrow g_{ij} \left( E_i(k) \otimes E_j(k) + F_i(k) \otimes F_j(k) \right)$$
$$b_{ij} \left( f_i e_j - e_i f_j \right) \rightarrow b_{ij} \left( F_i(k) \otimes E_j(k) - E_i(k) \otimes F_j(k) \right)$$

Finally, using the rule (7-i), the RHS of (20) is obtained.

### 2) DTs of the Designed Differential Equations

For the differential equation in **Formulation** 1, i.e., (14), its DT is derived as follows. After applying the rule (7-iv) for LHS and (7-v) for RHS, $(k+1)\Lambda(k+1)=C_1\delta(k)$ holds. Then, it can be further written in (24) after replacing $k$ by $k$-1 and using the definition of $\delta(k)$.

$$\Lambda(k) = C_1\delta(k-1) \quad (24)$$

For **Formulation** 2 in (15), the DT of the differential equation is in (25) where the derivation is similar as (24) and is omitted here; the DT of the ancillary equation is in (26) after applying the rule (7-iii) to both sides.

$$V_i(k) = -C_2\delta(k-1) \quad (25)$$

$$V_i(k) \otimes V_i(k) = E_i(k) \otimes E_i(k) + F_i(k) \otimes F_i(k) \quad (26)$$

### E. Step 3: Proving the Nonlinear Power Flow Equations Satisfy Formally Linear Equations after DT

**Proposition 1**: The transformed power flow equations (20)-(23) respectively satisfy formally linear equations (27)-(30).

$$0 = \boldsymbol{a}_{P,i} \boldsymbol{Y}(k) - \Delta p_i \Lambda(k) + \varepsilon_i, \text{if } i \in \Omega_{PQ} \cup \Omega_{PV} \quad (27)$$

$$0 = \boldsymbol{a}_{Q,i} \boldsymbol{Y}(k) - \Delta q_i \Lambda(k) + \mu_i, \text{if } i \in \Omega_{PQ} \quad (28)$$

$$0 = \boldsymbol{a}_{V,i} \boldsymbol{Y}(k) + 0\Lambda(k) + \varsigma_i, \text{if } i \in \Omega_{PV} \quad (29)$$

$$0 = \boldsymbol{a}_{E,i} \boldsymbol{Y}(k) + 0\Lambda(k) - e_i^{sp}\delta(k), \text{if } i \in \Omega_{REF}$$
$$0 = \boldsymbol{a}_{F,i} \boldsymbol{Y}(k) + 0\Lambda(k) - f_i^{sp}\delta(k), \text{if } i \in \Omega_{REF} \quad (30)$$

where $\boldsymbol{Y}(k) \in \mathbb{R}^{2N \times 1}$ and $\Lambda(k) \in \mathbb{R}$ are variables representing the DT of $\boldsymbol{y}$ and $\lambda$ respectively; $\boldsymbol{a}_{P,i}, \boldsymbol{a}_{Q,i}, \boldsymbol{a}_{V,i}, \boldsymbol{a}_{E,i}, \boldsymbol{a}_{F,i} \in \mathbb{R}^{1 \times 2N}$ and $\varepsilon_i$, $\mu_i$, $\zeta_i \in \mathbb{R}$ are parameters given in (44)-(49) respectively. The detailed proof of Proposition 1 is in Appendix.

From the **Proposition**, DTs (9c) of the nonlinear power flow equation satisfy formally linear equation (10a) with matrices $\boldsymbol{A}_{gy}$, $\boldsymbol{A}_{g\lambda}$, and $\boldsymbol{B}_g$ given by (31). For notation simplicity, here we let buses 1 to $M$ be PQ buses, buses $M$+1 to $N$-1 be PV buses and bus $N$ be the reference bus.

$$\boldsymbol{A}_{gy} = \begin{bmatrix} \boldsymbol{A}_{y,PQ} \\ \boldsymbol{A}_{y,PV} \\ \boldsymbol{A}_{y,REF} \end{bmatrix}, \boldsymbol{A}_{g\lambda} = \begin{bmatrix} \boldsymbol{A}_{\lambda,PQ} \\ \boldsymbol{A}_{\lambda,PV} \\ \boldsymbol{A}_{\lambda,REF} \end{bmatrix}, \boldsymbol{B}_g = \begin{bmatrix} \boldsymbol{B}_{PQ} \\ \boldsymbol{B}_{PV} \\ \boldsymbol{B}_{REF} \end{bmatrix} \quad (31)$$

$$\boldsymbol{A}_{y,PQ} = \begin{bmatrix} \boldsymbol{a}_{P,1} \\ \boldsymbol{a}_{Q,1} \\ \vdots \\ \boldsymbol{a}_{P,M} \\ \boldsymbol{a}_{Q,M} \end{bmatrix} \quad \boldsymbol{A}_{\lambda,PQ} = -\begin{bmatrix} \Delta p_1 \\ \Delta q_1 \\ \vdots \\ \Delta p_M \\ \Delta q_M \end{bmatrix} \quad \boldsymbol{B}_{PQ} = \begin{bmatrix} \varepsilon_1 \\ \mu_1 \\ \vdots \\ \varepsilon_M \\ \mu_M \end{bmatrix}$$

$$\boldsymbol{A}_{y,PV} = \begin{bmatrix} \boldsymbol{a}_{P,M+1} \\ \boldsymbol{a}_{V,M+1} \\ \vdots \\ \boldsymbol{a}_{P,N-1} \\ \boldsymbol{a}_{V,N-1} \end{bmatrix} \quad \boldsymbol{A}_{\lambda,PV} = -\begin{bmatrix} \Delta p_{M+1} \\ 0 \\ \vdots \\ \Delta p_{N-1} \\ 0 \end{bmatrix} \quad \boldsymbol{B}_{PV} = \begin{bmatrix} \varepsilon_{M+1} \\ \zeta_{M+1} \\ \vdots \\ \varepsilon_{N-1} \\ \zeta_{N-1} \end{bmatrix}$$

$$A_{y,\mathrm{REF}} = \begin{bmatrix} a_{\mathrm{E},N} \\ a_{\mathrm{F},N} \end{bmatrix} \quad A_{\lambda,\mathrm{REF}} = \begin{bmatrix} 0 \\ 0 \end{bmatrix} \quad B_{\mathrm{REF}} = \begin{bmatrix} -e_N^{sp}\delta(k) \\ -f_N^{sp}\delta(k) \end{bmatrix}$$

Besides, the DT (26) of the ancillary equation in Formulation II also satisfies a formally linear equation in (10b) with proof in the Appendix.

### F. Step 4: Non-iterative Algorithm for Solving Variables as Power Series of Time

Following the basic idea in Section III-B, an algorithm is designed to solve power series coefficients $X(k)$, $\Lambda(k)$, $Y(k)$ up to any desired order. Note that these coefficients are explicitly calculated with no numerical iteration.

---

**Algorithm: Solve Coefficients**

**Input**: initial values $x(0), y(0), \lambda(0)$

**Output**: any order coefficients $X(k), Y(k), \Lambda(k), k = 0 \cdots K$

**Steps:**
Initialization:
  $X(k = 0) = x(t = 0)$
  $Y(k = 0) = y(t = 0)$
  $\Lambda(k = 0) = \lambda(t = 0)$
  Calculate matrices $A_{gy}, A_{gz}, A_{hy}, A_{h\lambda}$ using (31)

  Calculate $A_{gy}^{-1}$ if Formulation I or $\begin{bmatrix} A_{hy} & A_{h\lambda} \\ A_{gy} & A_{g\lambda} \end{bmatrix}^{-1}$ if Formulation II

**for** $k = 1 : K$
  Calculate matrix $B_g, B_h$ using (31)
  **if** Formulation I
    Solve $\Lambda(k)$ using  $\Lambda(k) = C_1\delta(k-1)$
    Solve $Y(k)$ using  $Y(k) = -A_{gy}^{-1}(A_{g\lambda}\Lambda(k) + B_g)$
  **if** Formulation II
    Solve $V_t(k)$ using  $V_t(k) = -C_2\delta(k-1)$
    Solve $Y(k), \Lambda(k)$ using  $\begin{bmatrix} Y(k) \\ \Lambda(k) \end{bmatrix} = -\begin{bmatrix} A_{gy} & A_{g\lambda} \\ A_{hy} & A_{h\lambda} \end{bmatrix}^{-1}\begin{bmatrix} B_g \\ B_h \end{bmatrix}$
**end**

---

After the power series coefficients are calculated, $y(t)$ and $\lambda(t)$ are calculated by evaluating the power series of time in (13) and the solution curves are directly obtained. In practical, the multi-time window strategy [18]-[19] can be used to extend the convergence region of power series of time and ensure the accuracy. The time step length as well as the order $K$ of the power series of time are usually selected from trial simulations [19], and the impact of $K$ and time step length are also studied in [18]-[19].

The proposed method can also be applied to more complicated power system models such as 1) considering reactive power limit of generators, 2) ZIP load model. First, to consider the reactive power limit of generators, the proposed method can be slightly modified as follows: if a generator meets the reactive power limit, then it is changed from a PV bus to a PQ bus; correspondingly, the matrices $A$ and $B$ need to be re-calculated using (31). Later in Section IV-A, we demonstrate the proposed method for tracing PV curves considering reactive power limits. Second, for a nonlinear power flow equation with

ZIP loads, we proved in [29] that its DTs still satisfy formally linear equations. Therefore, the proposed method can be directly applied with slight modification on matrices $A$ and $B$.

Regarding the computational complexity, the proposed method has two unique features: First, it shifts most of the computation burden to the offline stage, i.e., deriving the equation for calculating matrices $A$ and $B$, which is a one-time effort (the matrices $A$ and $B$ derived in this paper can be directly used by others without deriving them again); and the online stage only involves explicit matrix operation and evaluation of analytical solutions, which do not require any numerical iterations. Second, the proposed method can reduce the frequency of solving linear equations compared with the CPF method, thus having better computational efficiency. This is because the CPF method needs to solve a linear equation in every iteration and every prediction-correction stage, while the proposed method only needs to solve linear equations once in each time step, and the total number of time steps are greatly reduced benefiting from the high order approximation.

## IV. CASE STUDIES

The proposed DPF method is first tested on the IEEE 9-bus system [13] to demonstrate the basic idea, the impact of load change directions, and the impact of reactive power limit of generators. Then, the accuracy, computational complexity, and computation time are compared with the CPF method in MATPOWER using several large systems including the IEEE 39-bus system, IEEE 57-bus system and a Polish 2383-bus test system [13]. At last, the proposed approach is applied to $N$-1 contingency analysis. Simulations are conducted in MATLAB R2017a on a personal computer with i5-8250U CPU. Without specification, generations and loads of all buses are uniformly increased. For the CPF method, various simulation control parameters are adjusted for the best time performance, including using the pseudo arc-length for parameterization, enabling adaptive step size, increasing the maximum allowed step size and disabling the incremental curve plotting in each iteration, etc. For the DPF method, parameters $C_1$ and $C_2$ are set as 1, $K$ is set as 6 from trail simulations, and the time step length is fixed at 0.05s for 2383-bus system and 0.1s for other systems.

### A. Demonstration on the 9-bus Power System

To demonstrate the idea of the proposed method, Fig. 2 and Fig. 3 respectively give the time domain trajectories of the solved dynamized power flow model and the obtained PV curve. In the first 1.63s, the loading parameter $\lambda$ increases with time in a constant rate while the voltage magnitude of bus 9 drops from 0.9956 p.u. to 0.6268 p.u., indicating high voltage solutions. During the time period between $t=1.63$s and $t=1.68$s, the voltage is decreased from 0.6268 p.u. to 0.5439 p.u., while the loading parameter is first increased from 1.63 to reach the maximum value 1.64 and then decreased to 1.63, indicating the dynamic process of passing the nose point. Finally, both the loading parameter and the bus voltage are decreased after $t=1.68$s, indicating the low voltage solutions. The obtained loading limit 1.64 is the same as the limit from the CPF method.

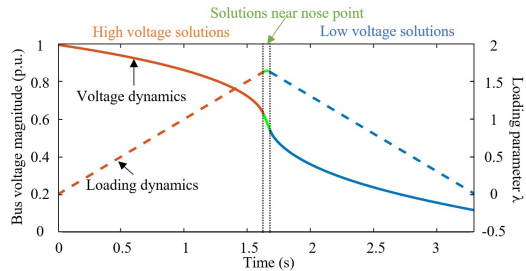

Fig. 2.  Time domain trajectory of the dynamized power flow model

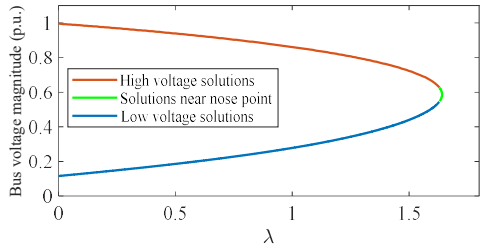

Fig. 3.  Solution curve of load bus 9 on 9-bus system

Two scenarios are designed to demonstrate the capability of the proposed method on handling load changes with 1) non-uniform directions and 2) reactive power limits. Fig. 4a shows the PV curve of load bus 9 when increasing generation at bus 3 and load at bus 7 by 50 MW in active power and 10 MVar in reactive power. Fig. 4b further shows the PV curve of the same bus when reactive power limit of generators is considered. It shows the calculated maximum loading limits are reduced from 8.17 to 7.79 due to the reactive power limit. These results demonstrate the performance of the proposed method on practical power system models and applications.

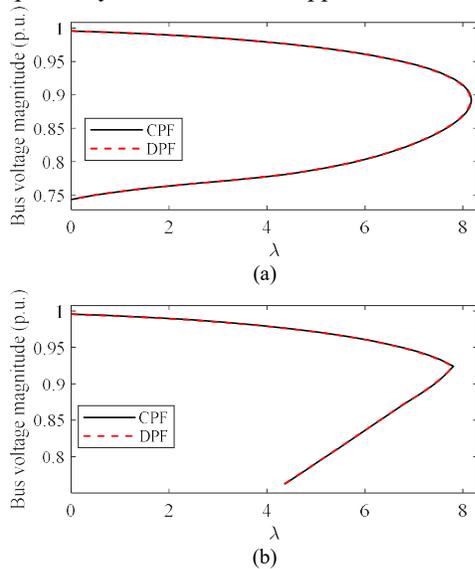

Fig. 4.  Solution curves under non-uniform load change direction, a) not consider reactive power limit, b) consider reactive power limit

### B.  Performance on Large Systems: Accuracy, Computation Complexity, and Time Performance

Respectively for the 39-bus system, the 57-bus system and the 2383-bus system, the proposed DPF method is compared with the CPF method. In all following studies, the CPF method is tested using the commercial MATPOWER package while the proposed DPF method is tested using our research code. Fig. 5 to Fig. 7 show the PV curves of three load buses, obtained by both the proposed method and the CPF method. Respectively for the three test systems, the calculated loading limits are 1.12, 0.88 and 0.89 for DPF method, and 1.13, 0.89 and 0.89 for the CPF method. These results demonstrate the accuracy of the proposed method.

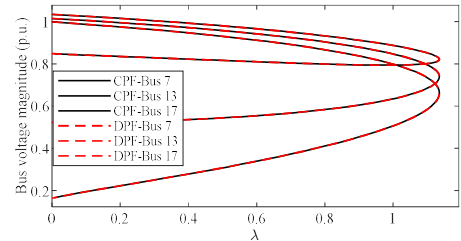

Fig. 5.  Solution curves on 39-bus system

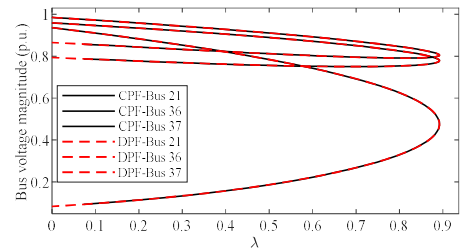

Fig. 6.  Solution curves on 57-bus system

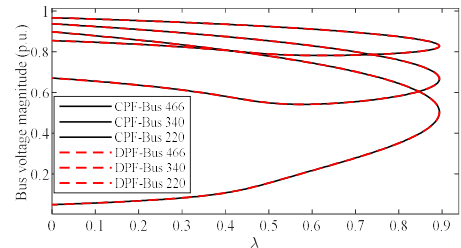

Fig. 7.  Solution curves on 2383-bus system

For both the CPF method and the DPF method, a major computation burden is in solving linear equations. Table I gives how many times linear equations are solved for both methods. It shows that the proposed approach is 10 times fewer than the CPF method for all the three test systems. This is because the CPF method needs to solve a linear equation in each iteration and for every prediction-correction step while the proposed method only solves a linear equation once in each time step.

TABLE I
NUMBERS OF TIMES OF SOLVING LINEAR EQUATIONS

| Test Systems | CPF | DPF |
|---|---|---|
| 39-bus system | 174 | 11 |
| 57-bus system | 108 | 10 |
| 2383-bus system | 424 | 18 |

Table II further gives the computation times of both methods. It shows the proposed DPF method is around 9 times, 12 times, and 2 times faster than the CPF method, respectively, for the three test systems. The speed up on the 2383-bus system is less than speedups on the other two smaller systems because our current academic research code that implements the DPF method in MATLAB has not been optimized to as efficiently handle large-scale matrix operations as the commercial CPF solver in the MATPOWER. However, these test results do demonstrate the potential of the proposed DPF method for online power flow solution tracing and voltage stability assessment.

TABLE II
COMPARISON OF TIME PERFORMANCE (UNIT: SECOND)

| Test Systems | CPF | DPF |
|---|---|---|
| 39-bus system | 0.26 | 0.03 |
| 57-bus system | 0.50 | 0.04 |
| 2383-bus system | 24.45 | 10.13 |

### C. Application to N-1 Contingency Analysis

The proposed approach is further applied to screen *N*-1 contingencies. For the 39-bus system, 46 contingencies are created each with the loss of each single branch. Fig. 8 shows the maximum loading condition identified by both methods. Using the CPF results as benchmarks, the DPF method is accurate and reliable for all the contingencies. Regarding the computation time, the CPF method and the DPF method respectively takes 12.0 s and 1.4 s, showing that the DPF method can identify insecure contingencies much faster than the CPF method, and thus can scan more contingencies than the CPF method within limited time in the real-time environment.

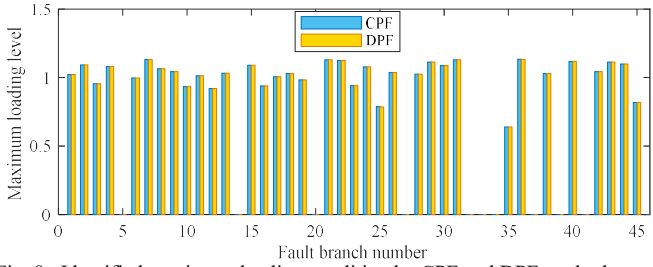

Fig. 8. Identified maximum loading condition by CPF and DPF method

## V. CONCLUSION

In this paper, a novel dynamized power flow method has been proposed to efficiently trace solution curves of power flow equations. The original curve tracing problem for steady-state power flow solutions is converted to a time domain simulation problem about a dynamized model after adding a differential equation on changes of the operating condition. An DT-based approach is proposed for efficiently solving the dynamized model without numerical iterations. Simulation results have shown high accuracy, reduced computational complexity and improved time performance of the proposed DPF method compared with a CPF solver in MATPOWER. Besides, the proposed method can deal with practical engineering constraints such as the non-uniform load change directions and reactive power limits of generators.

## APPENDIX

To make the proofs more compact, the following Lemma is first proved. In the Lemma, the transformation of multiplication operation from time domain to the convolution operation in the domain of power series orders is well-known in many DT literatures, however, the resulted linear relationship in (33)-(34), despite their simplicity and being straightforward, are rarely noticed and exploited as far as the authors know.

**Lemma:** The DT of $z(t)=x(t)y(t)$, satisfies a formally linear equation in (32). Especially, when $x(t)=y(t)$, (33) holds.

$$Z(k) = X(k) \otimes Y(k) = aX(k) + bY(k) + c \quad (32)$$

$$Z(k) = X(k) \otimes X(k) = 2aX(k) + c \quad (33)$$

**Proof of Lemma:**

$$Z(k) = X(k) \otimes Y(k) = \sum_{m=0}^{k} X(m)Y(k-m)$$

$$= X(0)Y(k) + X(k)Y(0) + \sum_{m=1}^{k-1} X(m)Y(k-m)$$

Therefore, (32) holds with *a*, *b* and *c* given below.

$$a = Y(0), b = X(0), c = \sum_{m=1}^{k-1} X(m)Y(k-m)$$

**Proof of Proposition 1:**

Use (27) as an example. The RHS of (20) is rewritten as

$$\text{RHS} = -\Delta p_i \Lambda(k) + \underbrace{g_{ii}\left(E_i(k) \otimes E_i(k) + F_i(k) \otimes F_i(k)\right)}_{\text{Term 1}}$$

$$+ \underbrace{\sum_{j=1, j\neq i}^{N} g_{ij}\left(E_i(k) \otimes E_j(k) + F_i(k) \otimes F_j(k)\right)}_{\text{Term 2}}$$

$$+ \underbrace{\sum_{j=1}^{N} b_{ij}\left(F_i(k) \otimes E_j(k) - E_i(k) \otimes F_j(k)\right)}_{\text{Term 3}}$$

According to the Lemma, the three terms are rewritten as:

$$\text{Term 1} = 2g_{ii}E_i(0)E_i(k) + 2g_{ii}F_i(0)F_i(k)$$

$$+ g_{ii}\sum_{m=1}^{k-1} E_i(m)E_i(k-m) + g_{ii}\sum_{m=1}^{k-1} F_i(m)F_i(k-m)$$

$$\text{Term 2} = \sum_{\substack{j=1, \\ j\neq i}}^{N} g_{ij}\left(E_j(0)E_i(k) + E_i(0)E_j(k)\right)$$

$$+ \sum_{\substack{j=1, \\ j\neq i}}^{N} g_{ij}\left(F_j(0)F_i(k) + F_i(0)F_j(k)\right)$$

$$+ \sum_{\substack{j=1, \\ j\neq i}}^{N} g_{ij}\left(\sum_{m=1}^{k-1} E_i(m)E_j(k-m) + \sum_{m=1}^{k-1} F_i(m)F_j(k-m)\right)$$

$$\text{Term 3} = \sum_{j=1}^{N} b_{ij}\left(E_j(0)F_i(k) + F_i(0)E_j(k)\right)$$

$$- \sum_{j=1}^{N} b_{ij}\left(F_j(0)E_i(k) + E_i(0)F_j(k)\right)$$

$$+ \sum_{j=1}^{N} b_{ij}\left(\sum_{m=1}^{k-1} F_i(m)E_j(k-m) - \sum_{m=1}^{k-1} E_i(m)F_j(k-m)\right)$$

Finally, (27) is obtained by summating the above three terms, with vector $\boldsymbol{a}_{P,i}$ and parameter $\varepsilon_i$ in (34) and (39). Similarly, (28)-(30) can be proved with vectors $\boldsymbol{a}_{P,i}$, $\boldsymbol{a}_{Q,i}$, $\boldsymbol{a}_{V,i}$, $\boldsymbol{a}_{E,i}$, $\boldsymbol{a}_{F,i}$ and parameters $\varepsilon_i$, $\mu_i$, $\zeta_i$ in (34)-(41).

$$\boldsymbol{a}_{\mathrm{P},i} = \begin{bmatrix} \alpha_{i1} & \beta_{i1} & \cdots & \alpha_{ij} & \beta_{ij} & \cdots \end{bmatrix}, \text{where}$$

$$\alpha_{ij} = g_{ij}E_i(0) + b_{ij}F_i(0), \ \ \beta_{ij} = g_{ij}F_i(0) - b_{ij}E_i(0), \text{if } j \neq i$$

$$\alpha_{ii} = \sum_{j=1}^{N}\left(g_{ij}E_j(0) - b_{ij}F_j(0)\right) + g_{ii}E_i(0) + b_{ii}F_i(0) \qquad (34)$$

$$\beta_{ii} = \sum_{j=1}^{N}\left(b_{ij}E_j(0) + g_{ij}F_j(0)\right) - b_{ii}E_i(0) + g_{ii}F_i(0)$$

$$\boldsymbol{a}_{\mathrm{Q},i} = \begin{bmatrix} \phi_{i1} & \psi_{i1} & \cdots & \phi_{ij} & \psi_{ij} & \cdots \end{bmatrix}, \text{where}$$

$$\phi_{ij} = -b_{ij}E_i(0) + g_{ij}F_i(0), \ \ \psi_{ij} = -b_{ij}F_i(0) - g_{ij}E_i(0), \text{if } j \neq i$$

$$\phi_{ii} = -\sum_{j=1}^{N}\left(b_{ij}E_j(0) + g_{ij}F_j(0)\right) - b_{ii}E_i(0) + g_{ii}F_i(0) \qquad (35)$$

$$\psi_{ii} = \sum_{j=1}^{N}\left(g_{ij}E_j(0) - b_{ij}F_j(0)\right) - g_{ii}E_i(0) - b_{ii}F_i(0)$$

$$\boldsymbol{a}_{\mathrm{V},i} = \begin{bmatrix} 0 & \cdots & 0 & 2E_i(0) & 2F_i(0) & 0 & \cdots & 0 \end{bmatrix} \qquad (36)$$

$$\boldsymbol{a}_{\mathrm{E},i} = \begin{bmatrix} 0 & \cdots & 0 & 1 & 0 & 0 & \cdots & 0 \end{bmatrix} \qquad (37)$$

$$\boldsymbol{a}_{\mathrm{F},i} = \begin{bmatrix} 0 & \cdots & 0 & 0 & 1 & 0 & \cdots & 0 \end{bmatrix} \qquad (38)$$

$$\varepsilon_i = \sum_{j=1}^{N} g_{ij}c_{ij} + \sum_{j=1}^{N} b_{ij}d_{ij} - p_i\delta(k), \text{where}$$

$$c_{ij} := \sum_{m=1}^{k-1} E_i(m)E_j(k-m) + \sum_{m=1}^{k-1} F_i(m)F_j(k-m) \quad (39)$$

$$d_{ij} := \sum_{m=1}^{k-1} F_i(m)E_j(k-m) - \sum_{m=1}^{k-1} E_i(m)F_j(k-m)$$

$$\mu_i = -\sum_{j=1}^{N} b_{ij}c_{ij} + \sum_{j=1}^{N} g_{ij}d_{ij} - q_i\delta(k) \qquad (40)$$

$$\varsigma_i = c_{ii} - v_i^2\delta(k) \qquad (41)$$

**Proof of (10b) from (26):** From the **Lemma**, there are

$$E_l(k) \otimes E_l(k) = 2E_l(0)E_l(k) + \sum_{m=1}^{k-1} E_l(m)E_l(k-m)$$

$$F_l(k) \otimes F_l(k) = 2F_l(0)F_l(k) + \sum_{m=1}^{k-1} F_l(m)F_l(k-m)$$

Then, (26) is rewritten as:

$$V_l(k) \otimes V_l(k) = 2E_l(0)E_l(k) + 2F_l(0)F_l(k)$$
$$+ \sum_{m=1}^{k-1} E_l(m)E_l(k-m) + \sum_{m=1}^{k-1} F_l(m)F_l(k-m)$$

Finally, (10b) is obtained with $\boldsymbol{a}_l$ and $\xi_l$ given in (42).

$$\boldsymbol{a}_l = \begin{bmatrix} 0 & \cdots & 0 & 2E_l(0) & 2F_l(0) & 0 & \cdots & 0 \end{bmatrix}$$
$$\xi_l = \sum_{m=1}^{k-1} E_l(m)E_l(k-m) + \sum_{m=1}^{k-1} F_l(m)F_l(k-m) - V_l(k) \otimes V_l(k) \qquad (42)$$